\documentclass{elsart}
\usepackage{amssymb}
\newtheorem{lemma}{Lemma}[section]
\newtheorem{theorem}{Theorem}[section]\newtheorem{proposition}{Proposition}[section]
\newtheorem{definition}{Definition}[section]

\begin{document}
\begin{frontmatter}
\title{On a mixed cubic-superlinear non radially symmetric Schr\"odinger system - Part I: Existence and uniqueness of solutions}
\author{Abdurahman F. Aljohani}
\ead{a.f.aljohani@ut.edu.sa}
%\corauth[cor1]{Corresponding author}
\address{Department of Mathematics, Faculty of Science, University of Tabuk, Saudi Arabia.}
%%%%%%%%%%%%%%%%
\author{Anouar Ben Mabrouk\corauthref{cor1}}
\ead{anouar.benmabrouk@fsm.rnu.tn; amabrouk@ut.edu.sa}
\corauth[cor1]{Corresponding author}
\address{Research Unit of Algebra, Number Theory and Nonlinear Analysis UR11ES50, Department of Mathematics, Faculty of Sciences, 5019 Monastir. Tunisia.\\
\&
Department of Mathematics, Higher Institute of Applied Mathematics and Informatics, Street of Assad Ibn Al-Fourat, Kairouan 3100, Tunisia.\\
\& Department of Mathematics, Faculty of Science, University of Tabuk, Saudi Arabia.}
%%%%%%%%%%%%%%%%%%%%%%%%%%%%%%%%%%%%%%%%%%%%%%%%%%%%%%%%%%%%%%%%%%%%%%%%%%%%%%%%%%%%%%%%%%%%%%%%%%%%%%%%%%%%%%%%%%%%%%%%%
%%%%%%%%%%%%%%%%%%%%%%%%%%%%%%%%%%%%%%%%%%%%%%%%%%%%%%%%%%%%%%%%%%%%%%%%%%%%%%%%%%%%%%%%%%%%%%%%%%%%%%%%%%%%%%%%%%%%%%%%%
%%%%%%%%%%%%%%%%%%%%%%%%%%%%%%%%%%%%%%%%%%%%%%%%%%%%%%%%%%%%%%%%%%%%%%%%%%%%%%%%%%%%%%%%%%%%%%%%%%%%%%%%%%%%%%%%%%%%%%%%%
\begin{abstract}
In this paper a nonlinear coupled Schrodinger system in the presnce of mixed cubic and superlinear power laws is considered. Focus are made on the steady state solutions of the continuous system for existence and uniqueness by minimizing problem on some associated Nehari manifold. 
\end{abstract}
\end{frontmatter}
\section{Introduction}
The present work is devoted to the study of a coupled system of nonlinear Schr\"odinger equations characterized by a mixed nonlinearities. Focuses are made on a theoretical study concerned with existence and uniqueness of solutions especially positive ones of the steady state problem associated to the evolutive one. 

Denote for $\lambda$ and $p$ real numbers such that $\lambda>0$ and $p>1$,
$$
g(u,v)=|u|^{p-1}+\lambda|v|^{2}\;\mbox{and}\;f(u,v)=g(u,v)u.
$$
We consider in the first part the evolutive system
\begin{equation}\label{ContinuousProblem1}
\left\{\begin{array}{lll}
iu_{t}+\sigma_1\Delta u+g(u,v)u=0,\\
iv_{t}+\sigma_2\Delta v+g(v,u)v=0
\end{array}
\right.
\end{equation}
on suitable domain $\mathcal{Q}=\Omega\times(t_0,+\infty)$, where $\Omega$ is a domain in $\mathbb{R}^N$, $N\geq1$ and $t_0$ is a reel parameter considered as the initial time. $u_{t}$ is the first order partial derivative in time, $\Delta=\displaystyle\sum_{i=1}^{N}\frac{\partial ^2}{\partial x_i^2}$ is the Laplace operator on $\mathbb{R}^N$. $\sigma_i$, $i=1,2$ are real parameters such that $\sigma_i>0$.

In the next section (section 2), a full study of existence and uniqueness of solutions of the steady state problem associated to problem (\ref{ContinuousProblem1}) on the whole space $\Omega=\mathbb{R}^N$ is developed. The main tool is by applying the notion of Nehari manifolds.
\section{The steady state problem}
Recall that a steady state solution of problem (\ref{ContinuousProblem1}) is any solution of the form
$$
W(x,y,t)=(e^{i\omega_1t}u(x,y),e^{i\omega_2t}v(x,y)).
$$
In a first step we claim that whenever $W$ is a steady state solution of problem (\ref{ContinuousProblem1}), then $(u,v)$ is a solution of the system
\begin{equation}\label{ContinuousSteadyStateProblem1}
\left\{\begin{array}{lll}
-\omega_1u+\sigma_1\Delta u+g(u,v)u=0,\\
-\omega_2v+\sigma_2\Delta v+g(v,u)v=0.
\end{array}
\right.
\end{equation}
This follows easily by substituting $W$ in the system (\ref{ContinuousProblem1}).

In the rest of the paper we will assume for simplicity that $\omega_1=\omega_2=\omega>0$. Therefore, $(u,v)$ is a solution of the system
\begin{equation}\label{ContinuousSteadyStateProblem2}
\left\{\begin{array}{lll}
-\omega u+\sigma_1\Delta u+g(u,v)u=0,\\
-\omega v+\sigma_2\Delta v+g(v,u)v=0.
\end{array}
\right.
\end{equation}
To study existence of the solutions of the last problem, we will apply the so-called Nehari manifolds. Consider the functional space
$$
H=H^1(\mathbb{R}^N)\times H^1(\mathbb{R}^N)
$$ equipped with the inner product
$$
<(u,v);(\varphi,\psi)>_H=\displaystyle\int_{\mathbb{R}^N}(\sigma_1\nabla u\nabla\varphi+\omega u\varphi+\sigma_2\nabla v\nabla\psi+\omega v\psi)dxdy
$$
and the norm
$$
\|(u,v)\|_H=\Bigl(<(u,v);(u,v)>_H\Bigr)^{1/2}.
$$
\begin{lemma}
$\Bigl(H,\|(u,v)\|_H\Bigr)$ is a Hilbert space.
\end{lemma}
\begin{definition}
A pair $w=(u,v)$ is said to be a weak solution of problem (\ref{ContinuousSteadyStateProblem2}) if $w\in H$ and for all $(\varphi,\psi)\in\mathcal{C}^{\infty}_c(\mathbb{R}^2)\times\mathcal{C}^{\infty}_c(\mathbb{R}^2)$,
$$
\displaystyle\int_{\mathbb{R}^2}(\sigma_1\nabla u\nabla\varphi+\omega u\varphi+\sigma_2\nabla v\nabla\psi+\omega v\psi-f(u,v)\varphi-f(v,u)\psi)dxdy=0
$$
where
$$
f(u,v)=(|u|^{p-1}+\lambda|v|^{2})u.
$$
\end{definition}
It is obvious that $(0,0)$ is a solution of (\ref{ContinuousSteadyStateProblem2}). We are interested in findiing nontrivial solutions for problem (\ref{ContinuousSteadyStateProblem2}) by applying critical point theory and the Nehari manifold. Consider the functional
$$
I(u,v)=\displaystyle\frac{1}{2}\mathcal{L}(u,v)-
\displaystyle\frac{1}{p+1}\mathcal{M}_p(u,v)-\displaystyle\frac{1}{2}\mathcal{N}_\lambda(u,v),
$$
where
$$
\mathcal{L}(u,v)=\displaystyle\int_{\mathbb{R}^N}\Bigl(\sigma_1|\nabla u|^2+\sigma_2|\nabla v|^2+\omega(u^2+v^2)\Bigr)dx,
$$
$$
\mathcal{M}_p(u,v)=\displaystyle\int_{\mathbb{R}^N}\Bigl(|u|^{p+1}+|v|^{p+1}\Bigr)dx
$$
and
$$
\mathcal{N}_\lambda(u,v)=\lambda\displaystyle\int_{\mathbb{R}^N}|u|^2|v|^2dx.
$$
The following holds.
\begin{lemma}
	\begin{description}
		\item[(i.)] The functional $I$ is $\mathcal{C}^1(H,\mathbb{R})$.
		\item[(ii.)] Any critical point of $I$ on $H$ is a weak solution pair of (\ref{ContinuousSteadyStateProblem2}).
		\item[(iii.)] The functional $I$ is not bounded neither above nor below on $H$.
	\end{description}
\end{lemma}
\textbf{Proof.}
	Straightforward calculus yield that
	$$
	\mathcal{L}(u+\varphi,v+\psi)=\mathcal{L}(u,v)+2\displaystyle\int_{\mathbb{R}^N}\Bigl(\sigma_1\nabla u\nabla\varphi+\sigma_2\nabla v\nabla\psi+\omega(u\varphi+v\psi)\Bigr)dx+....
	$$
	Hence, $\mathcal{L}$ is differentiable and
	$$
	\mathcal{L}'(u,v)(\varphi,\psi)=2\displaystyle\int_{\mathbb{R}^N}\Bigl(\sigma_1\nabla u\nabla\varphi+\sigma_2\nabla v\nabla\psi+\omega(u\varphi+v\psi)\Bigr)dx.
	$$
	Next,
	$$
	\mathcal{M}_p(u+\varphi,v+\psi)=\mathcal{M}_p(u,v)+(p+1)\displaystyle\int_{\mathbb{R}^N}\Bigl(|u|^{p-1}u\varphi+|v|^{p-1}v\psi\Bigr)dx+....
	$$
	Consequently, $\mathcal{M}_p$ is differentiable and
	$$
	\mathcal{M}_p'(u,v)(\varphi,\psi)=(p+1)\displaystyle\int_{\mathbb{R}^N}\Bigl(|u|^{p-1}u\varphi+|v|^{p-1}v\psi\Bigr)dx.
	$$
	Finally and similarly we get
	$$
	\mathcal{N}_\lambda(u+\varphi,v+\psi)=\mathcal{N}_\lambda(u,v)+2\lambda\displaystyle\int_{\mathbb{R}^N}(|u||v|^2\varphi+|v||u|^2\psi)dx+....
	$$
	Hence, $\mathcal{N}_\lambda$ is differentiable and
	$$
	\mathcal{N}_\lambda'(u,v)(\varphi,\psi)=2\lambda\displaystyle\int_{\mathbb{R}^N}(|u||v|^2\varphi+|v||u|^2\psi)dx.
	$$
	As a result, $I$ is differentiable and
	$$
	\begin{array}{lll}
	I'(u,v)(\varphi,\psi)&=&\mathcal{L}'(u,v)(\varphi,\psi)-\mathcal{M}_p'(u,v)(\varphi,\psi)-\mathcal{N}_\lambda'(u,v)(\varphi,\psi)\\
	&=&\displaystyle\int_{\mathbb{R}^N}\Bigl(\sigma_1\nabla u\nabla\varphi+\sigma_2\nabla v\nabla\psi+\omega(u\varphi+v\psi)\Bigr)dx\\
	&&-\displaystyle\int_{\mathbb{R}^N}\Bigl(|u|^{p-1}u\varphi+|v|^{p-1}v\psi\Bigr)dx\\
	&&-\lambda\displaystyle\int_{\mathbb{R}^N}(|u||v|^2\varphi+|v||u|^2\psi)dx.
	\end{array}
	$$

To overcome the problem of the boundless of the functional $I$, we will apply instead the Nehari manifold notion to prove the existence of at least one positive solution of problem (\ref{ContinuousSteadyStateProblem2}). To do this, we consider the Nehari manifold
$$
\mathcal{NH}=\left\{(u,v)\in H\setminus\{0\};\biggl<I'(u,v),(u,v)\biggr>=0\right\}.
$$
The following result holds.
\begin{proposition}
	The Nehari manifold $\mathcal{NH}$ is not empty.
\end{proposition}
\textbf{Proof.}
	Let next $(u,v)\in H$ such that $u>0$ and $v>0$ a.e. on $\mathbb{R}^N$ and consider for $t\geq0$,
	$$
	H_{u,v}(t)=I(tu,tv)=\displaystyle\frac{t^2}{2}\mathcal{L}(u,v)-\displaystyle\frac{t^{p+1}}{p+1}\mathcal{M}_p(u,v)-\displaystyle\frac{t^4}{2}\mathcal{N}_\lambda(u,v).
	$$
	It is straightforward that $H_{u,v}\in\mathcal{C}^{2}(\mathcal{R}_+^*,\mathcal{R})$ and
	$$
	H'_{u,v}(t)=t\Bigl(\mathcal{L}(u,v)-t^{p-1}\mathcal{M}_p(u,v)-2t^2\mathcal{N}_\lambda(u,v)\Bigr).
	$$
	Denote
	$$
	\mathcal{K}_{u,v}(t)=\mathcal{L}(u,v)-t^{p-1}\mathcal{M}_p(u,v)-2t^2\mathcal{N}_\lambda(u,v).
	$$
	It is easy to see that $\mathcal{K}_{u,v}$ is a strictly decreasing function on $(0,+\infty)$. Furthermore,
	$$
	\mathcal{K}_{u,v}(0)=\mathcal{L}_{u,v}>0\quad\hbox{and}\quad\displaystyle\lim_{t\rightarrow+\infty}\mathcal{K}_{u,v}(t)=-\infty.
	$$
	Next, there exists $t_0>0$ for which
	$$
	\mathcal{K}_{u,v}(t_0)=0.
	$$
	At this point $t_0$, we get
	$$
	\mathcal{H}'_{u,v}(t_0)=0.
	$$
	On the other hand,
	$$
	\mathcal{H}'_{u,v}(t_0)=\Bigl<I'(t_0u,t_0v),(u,v)\Bigr>.
	$$
	Hence,
	$$
	\Bigl<I'(t_0u,t_0v),(u,v)\Bigr>=0.
	$$
	Denote next $\varphi=t_0u$ and $\psi=t_0v$. It is obvious that $(\varphi,\psi)\in H\setminus\{0\}$. Else, we have
	$$
	\Bigl<I'(\varphi,\psi),(\varphi,\psi)\Bigr>=0.
	$$
	Consequently, $(\varphi,\psi)\in\mathcal{NH}$.

\begin{proposition}\label{boundedseq}
	For any sequence $(u_n,v_n)_n$ in $\mathcal{NH}$, the following assertion holds.
	\begin{quote}
		$(I(u_n,v_n))_n$ is bounded (in $\mathbb{R}$) if and only if $(u_n,v_n)_n$ is bounded in $H$.
	\end{quote}
\end{proposition}
\textbf{Proof.} We claim that
	\begin{equation}\label{Iunvnclaim1}
	I(u_n,v_n)\geq\max(\displaystyle\frac{p-1}{2(p+1)},\displaystyle\frac{1}{4})\|(u_n,v_n)\|_H^2;\;\;\forall(u_n,v_n)_m\subset\mathcal{NH}.
	\end{equation}
	Indeed, recall that
	\begin{equation}\label{Iunvn1}
	I(u_n,v_n)=\displaystyle\frac{1}{2}\mathcal{L}(u_n,v_n)-
	\displaystyle\frac{1}{p+1}\mathcal{M}_p(u_n,v_n)-\displaystyle\frac{1}{2}\mathcal{N}_\lambda(u_n,v_n),
	\end{equation}
	On the other hand, $(u_n,v_n)_n\subset\mathcal{NH}$. Henceforth,
	$$
	I'(u_n,v_n)(u_n,v_n)=\mathcal{L}(u_n,v_n)-\mathcal{M}_p(u_n,v_n)-2\mathcal{N}_\lambda(u_n,v_n)=0.
	$$
	Therefore,
	\begin{equation}\label{Iunvn2}
	\mathcal{M}_p(u_n,v_n)=\mathcal{L}(u_n,v_n)-2\mathcal{N}_\lambda(u_n,v_n).
	\end{equation}
	Combining equations (\ref{Iunvn1}) and (\ref{Iunvn2}) we get
	\begin{equation}\label{Iunvn3}
	I(u_n,v_n)=\displaystyle\frac{p-1}{2(p+1)}\mathcal{L}(u_n,v_n)+\displaystyle\frac{3-p}{2(p+1)}\mathcal{N}_\lambda(u_n,v_n).
	\end{equation}
	As a result, whenever $p\leq3$, we get
	\begin{equation}\label{Iunvn4}
	I(u_n,v_n)\geq\displaystyle\frac{p-1}{2(p+1)}\mathcal{L}(u_n,v_n).
	\end{equation}
	Now, for $p>3$, by coming back to (\ref{Iunvn2}) we get
	\begin{equation}\label{Iunvn5}
	\mathcal{L}(u_n,v_n)=\mathcal{M}_p(u_n,v_n)+2\mathcal{N}_\lambda(u_n,v_n).
	\end{equation}
	Hence,
	$$
	\mathcal{L}(u_n,v_n)\geq2\mathcal{N}_\lambda(u_n,v_n).
	$$
	Which yields that
	\begin{equation}\label{Iunvn6}
	\displaystyle\frac{3-p}{2(p+1)}\mathcal{N}_\lambda(u_n,v_n)\geq\displaystyle\frac{3-p}{4(p+1)}\mathcal{L}(u_n,v_n).
	\end{equation}
	Now and finally, by combining equations (\ref{Iunvn3}) and (\ref{Iunvn6}) we get
	\begin{equation}\label{Iunvn7}
	I(u_n,v_n)\geq\Biggl(\displaystyle\frac{p-1}{2(p+1)}+\displaystyle\frac{3-p}{4(p+1)}\Biggr)\mathcal{L}(u_n,v_n)=\displaystyle\frac{1}{4}\mathcal{L}(u_n,v_n).
	\end{equation}
	The equations (\ref{Iunvn4}) and (\ref{Iunvn7}) yields that
	\begin{equation}\label{Iunvnclaim2}
	I(u_n,v_n)\geq\max(\displaystyle\frac{p-1}{2(p+1)},\displaystyle\frac{1}{4})\mathcal{L}(u_n,v_n);\;\;\forall(u_n,v_n)_m\subset\mathcal{NH}.
	\end{equation}
	Observing that
	$$
	\mathcal{L}(u_n,v_n)=\|(u_n,v_n)\|_H^2;\;\;\forall(u_n,v_n)_m\subset\mathcal{NH},
	$$
	claim (\ref{Iunvnclaim1}) holds. Hence, the necessary condition follows. We next show the sufficient one. So assume that the sequence $(u_n,v_N)_n\subset\mathcal{NH}$ is bounded in $H$. Equation (\ref{Iunvn5}) written otherwise yields that
	\begin{equation}\label{Iunvn8}
	\|(u_n,v_n)\|_H^2=\mathcal{M}_p(u_n,v_n)+2\mathcal{N}_\lambda(u_n,v_n).
	\end{equation}
	Consequently, both $(\mathcal{M}_p(u_n,v_n))_n$ and $(\mathcal{N}_\lambda(u_n,v_n))_n$ are bounded. Next, again equation (\ref{Iunvn1}) written otherwise yields that
	\begin{equation}\label{Iunvn9}
	I(u_n,v_n)=\displaystyle\frac{1}{2}\|(u_n,v_n)\|_H^2-\displaystyle\frac{1}{p+1}\mathcal{M}_p(u_n,v_n)-\displaystyle\frac{1}{2}\mathcal{N}_\lambda(u_n,v_n).
	\end{equation}
	So $(I(u_n,v_n))_n$ is a linear combination of bounded sequences. Hence, it is also bounded.

We now state our first result on the existence of solutions for problem (\ref{ContinuousSteadyStateProblem2}).
\begin{theorem}
	Problem (\ref{ContinuousSteadyStateProblem2}) has a nontrivial positive solution.
\end{theorem}
\textbf{Proof.}
	We proceed by steps. We claim the following.\\
	\textit{Claim 1.} $I_0=\displaystyle\inf\{I(u,v)\,;\;\;(u,v)\in\mathcal{NH}\}>0$.\\
	\textit{Claim 2.} There exists sequences $(u_n,v_n)_n\subset\mathcal{NH}$ and $(a_n)_n\subset\mathbb{R}$ such that
	\begin{equation}\label{claim2sequence}
	\left\{\begin{array}{lll}
	I(u_n,v_n){\longrightarrow}I_0\\
	\mbox{and}\\
	I'(u_n,v_n)-a_nF'(u_n,v_n){\longrightarrow}0\mbox{ in }H',
	\end{array}\right.
	\end{equation}
	as ${n\to+\infty}$, where $H'$ is the dual of $H$.\\
	\textit{Claim 3.} Any sequence satisfying (\ref{claim2sequence}) is a Palais-Smale sequence of $I$ in $H$.\\
	\textit{Claim 4.} There exists a sequence $(u_n,v_n)_n\subset\mathcal{NH}$ that is a Palais-Smale sequence of $I$ in $H$.\\
	\textit{Claim 5.} There exists a pair $(u,v)\in H$ satisfying
	\begin{equation}\label{claim5}
	I(u,v)=I_0\quad\mbox{and}\quad I'(u,v)=0.
	\end{equation}
	\textit{Claim 6.} The pair $(|u|,|v|)$ is a solution of (\ref{ContinuousSteadyStateProblem2}), which is positive and where $(u,v)$ is defined in claim 5.

We now proceed in proving claims 1 to 6.\\
\textit{Proof of Claim 1.} Recall that (equation (\ref{Iunvnclaim1}))
\begin{equation}\label{Iunvnclaim1bis}
I(u_n,v_n)\geq\,C_p\|(u_n,v_n)\|_H^2;\;\;\forall(u_n,v_n)_n\subset\mathcal{NH},
\end{equation}
where
$$
C_p=\max(\displaystyle\frac{p-1}{2(p+1)},\displaystyle\frac{1}{4}).
$$
Therefore $I_0\geq0$. Assume next that $I_0=0$. There exists then a sequence $(u_n,v_n)_n\subset\mathcal{NH}$ such that
$$
I(u_n,v_n){\longrightarrow}0\;\mbox{as}\;{n\to+\infty}.
$$
Using again equation (\ref{Iunvnclaim1bis}), it follows that
$$
x_n=\|(u_n,v_n)\|_H{\longrightarrow}0\;\mbox{as}\;{n\to+\infty}.
$$
Next, from (\ref{Iunvn5}) it follows that
\begin{equation}\label{Iunvn5bis}
x_n^2=\mathcal{M}_p(u_n,v_n)+2\mathcal{N}_\lambda(u_n,v_n).
\end{equation}
It follows from Sobolev embedding theorem, there exists a constant $C_{p,N}>0$ such that
\begin{equation}\label{Iunvn5bis1}
\mathcal{M}_p(u_n,v_n)\leq\,C_{p,N}\|(u_n,v_n)\|_H^{p+1}.
\end{equation}
Similarly, there exists a constant $C_{\lambda,N}>0$ such that
\begin{equation}\label{Iunvn5bis2}
\mathcal{N}_\lambda(u_n,v_n)\leq\,C_{\lambda,N}\|(u_n,v_n)\|_H^{4}.
\end{equation}
As a result, by combining equations (\ref{Iunvn5bis}), (\ref{Iunvn5bis1}) and (\ref{Iunvn5bis2}), we get
\begin{equation}\label{Iunvn5bis3}
x_n^2\leq\,C_{p,N}x_n^{p+1}+C_{\lambda,N}x_n^{4}.
\end{equation}
Or equivalently,
\begin{equation}\label{Iunvn5bis3}
1\leq\,C_{p,N}x_n^{p-1}+C_{\lambda,N}x_n^{2}.
\end{equation}
Which contradicts the fact that $x_n{\longrightarrow}0$ as ${n\to+\infty}$.\\
\textit{Proof of Claim 2.} By definition of $I_0$, we have
$$
I_0=\inf\left\{I(u,v)\in H\setminus\{0\};F(u,v)=0\right\},
$$
where
$$
F(u,v)=\biggl<I'(u,v),(u,v)\biggr>.
$$
Hence, by applying Ekeland's variational principle, it follows that there exists $(u_n,v_n)_n\subset\mathcal{NH}$ satisfying (\ref{claim2sequence}). Hence, claim 2 holds.\\
\textit{Proof of Claim 3.} Let $(u_n,v_n)_n$ be satisfying claim 2. It suffices to show that
\begin{equation}\label{claim2sequence1}
a_nF'(u_n,v_n){\longrightarrow}0\mbox{ in }H'\;\mbox{as}\;{n\to+\infty}.
\end{equation}
Indeed, as
\begin{equation}
I(u_n,v_n){\longrightarrow}I_0\;\mbox{as}\;{n\to+\infty},
\end{equation}
we get by Proposition \ref{boundedseq} that $(u_n,v_n)_n$ is bounded in $H$. Hence,
\begin{equation}\label{boundedseq1}
\biggl<I'(u_n,v_n),(u_n,v_n)\biggr>-a_n\biggl<F'(u_n,v_n),(u_n,v_n)\biggr>{\longrightarrow}0\;\mbox{as}\;{n\to+\infty}.
\end{equation}
Oberving that
\begin{equation}\label{boundedseq2}
\biggl<I'(u_n,v_n),(u_n,v_n)\biggr>=0,
\end{equation}
it follows that
\begin{equation}\label{boundedseq3}
a_n\biggl<F'(u_n,v_n),(u_n,v_n)\biggr>{\longrightarrow}0\;\mbox{as}\;{n\to+\infty}.
\end{equation}
On the other hand,
\begin{equation}\label{boundedseq4}
\biggl<F'(u_n,v_n),(u_n,v_n)\biggr>=(1-p)\mathcal{M}_p(u_n,v_n).
\end{equation}
So, it is bounded due to (\ref{Iunvn8}). Consequently, equation (\ref{boundedseq3}) yields that
\begin{equation}\label{boundedseq5}
a_n{\longrightarrow}0\;\mbox{as}\;{n\to+\infty}.
\end{equation}
Henceforth, (\ref{claim2sequence}) implies that
\begin{equation}\label{boundedseq6}
I'(u_n,v_n){\longrightarrow}0\mbox{ in }H'\;\mbox{as}\;{n\to+\infty}.
\end{equation}
\textit{Proof of Claim 4.} It follows from claim 2 and claim 3 immediately.\\
\textit{Proof of Claim 5.} We will prove that the sequence $(u_n,v_n)_n$ of claim 4 is convergent in $H$ to a pair $(u,v)$ (up to a subsequence) satisfying (\ref{claim5}). Indeed, from Proposition \ref{boundedseq} the sequence $(u_n,v_n)_n$ is bounded in $H$. Hence, there exists a pair $(u,v)\in H$ such that
\begin{equation}\label{weakconvergenceunvn}
u_n{\longrightarrow}u\quad\mbox{and}\quad
v_n{\longrightarrow}v\;\mbox{as}\;{n\to+\infty}
\end{equation}
in $H^1(\mathbb{R}^N)$. We shall show that
\begin{equation}\label{boundedseq7}
I(u,v)=I_0\quad\mbox{and}\quad I'(u,v)=0.
\end{equation}
Let $(\varphi,\psi)\in\mathcal{C}_c^\infty(\mathbb{R}^N)$ be a test pair. It follows from (\ref{weakconvergenceunvn}) that
\begin{equation}\label{weakconvergenceunvn1}
\mathcal{L}'(u_n,v_n)(\varphi,\psi){\longrightarrow}\mathcal{L}'(u,v)(\varphi,\psi)\;\mbox{as}\;{n\to+\infty}.
\end{equation}
Next, as $(\varphi,\psi)\in\mathcal{C}_c^\infty(\mathbb{R}^N)$, it results that for ${n\to+\infty}$,
\begin{equation}\label{weakconvergenceunvn2}
\mathcal{M}_p'(u_n,v_n)(\varphi,\psi){\longrightarrow}\mathcal{M}_p'(u,v)(\varphi,\psi)
\end{equation}
and similarly, whenever ${n\to+\infty}$, it holds that 
\begin{equation}\label{weakconvergenceunvn3}
\mathcal{N}_\lambda'(u_n,v_n)(\varphi,\psi){\longrightarrow}\mathcal{N}_\lambda'(u,v)(\varphi,\psi).
\end{equation}
As a result, as ${n\to+\infty}$ we get
\begin{equation}\label{weakconvergenceunvn4}
0=\biggl<I'(u_n,v_n),(\varphi,\psi)\biggr>{\longrightarrow}\biggl<I'(u,v),(\varphi,\psi)\biggr>.
\end{equation}
Hence, the second part of (\ref{boundedseq7}). It remains to prove its first part. We already now that up to a subsequence, we have
\begin{equation}
\|(u_n,v_n)\|_H{\longrightarrow}\|(u,v)\|_H\;\mbox{as}\;{n\to+\infty}.
\end{equation}
So, as the functional $I$ is continuus we get the desired result.\\
\textit{Proof of Claim 6.} It follows from the parity of the functionals $I$ and $F$ that
\begin{equation}\label{|u||v|solution1}
I(|u|,|v|)=I_0\quad\mbox{and}\quad F(|u|,|v|)=0.
\end{equation}
Which means that $(|u|,|v|)\in\mathcal{NH}$. Therefore,
\begin{equation}\label{|u||v|solution2}
\exists\,a\in\mathbb{R};\;\;I'(|u|,|v|)=aF'(|u|,|v|).
\end{equation}
IF we succeed in proving that $a=0$, we get immediately $I'(|u|,|v|)=0$, which means that $(|u|,|v|)$ is a solution of problem (\ref{ContinuousSteadyStateProblem2}). Of course, it is positive. Now, similar techniques as for (\ref{Iunvn2}) and (\ref{boundedseq4}), we get
\begin{equation}\label{|u||v|solution3}
\bigl<F'(|u|,|v|),(|u|,|v|)\bigr>=(1-p)\mathcal{M}_p(|u|,|v|).
\end{equation}
Consequently, equation (\ref{|u||v|solution3}) yields that
\begin{equation}\label{|u||v|solution4}
0=\bigl<I'(|u|,|v|),(|u|,|v|)\bigr>=a(1-p)\mathcal{M}_p(|u|,|v|).
\end{equation}
Hence, $a=0$.
\section{Conclusion}
The present study of existence an uniqueness of positive solutions for a steady state nonlinear Schr\"odinger system is developed. The problem is considered in the presence of mixed nonlinearities such as a cubic power law and a superlinear one. The main tools are based on the notion of Nehari manifold and some variational techniques. We intend in a forthcoming part to develop a non standard numerical method to study numerical solutions of the original evolutive system by means of sophisticated algebraic operators such as the famous Lyapunov-Sylvester ones in a higher dimensional case. 

\end{document}